\documentclass[12pt,reqno]{amsart}

\usepackage{amssymb}
\usepackage{amsmath}
\usepackage{amsfonts}

\def\Z{\mathbb{Z}}

\def\F{\mathbb{F}}
\def\G{\mathcal{G}}
\def\P{\mathcal{P}}
\def\L{\mathcal{L}}
\newcommand{\Aut}{\mathrm{Aut}}

\begin{document}

\title{A new partial geometry $pg(5,5,2)$}

\author{Vedran Kr\v{c}adinac}

\address{Department of Mathematics, Faculty of Science, University of Zagreb,
Bijeni\v{c}ka~30, HR-10000 Zagreb, Croatia}

\email{vedran.krcadinac@math.hr}

\thanks{This work has been supported by the Croatian Science
Foundation under the projects $6732$ and $9752$.}



\begin{abstract}
We construct a new partial geometry with parameters $pg(5,5,2)$, not
isomorphic to the partial geometry of van Lint and Schrijver.
\end{abstract}

\maketitle

\section{Introduction}

A \emph{partial geometry} $pg(s,t,\alpha)$ is a partial linear space
with lines of degree $s+1$ and points of degree $t+1$ such that for
every non-incident point-line pair $(P,\ell)$, there are exactly
$\alpha$ points on $\ell$ collinear with $P$. Consequently, the
number of points is $v=(s+1)(st/\alpha+1)$ and the point graph is
strongly regular with parameters
$srg(v,s(t+1),s-1+t(\alpha-1),\alpha(t+1))$. The dual of a
$pg(s,t,\alpha)$ is a $pg(t,s,\alpha)$, hence formulae for the
number of lines and parameters of the line graph are obtained by
exchanging $s$ and $t$. For results about partial geometries we
refer to~\cite{CM95, DC03, JT07}.

Van Lint and Schrijver~\cite{vLS81} constructed a partial geometry
with parameters $pg(5,5,2)$. Another construction of the same
partial geometry was given in~\cite{CvL82}. This is one of only
three known proper partial geometries with $\alpha=2$, all of them
sporadic examples~\cite{DC03}. It was denoted as a partial geometry
of Type 6 in~\cite{JT07} and was characterized in~\cite{DW06} as the
only $pg(s,t,2)$ with an abelian Singer group of rigid type (i.e.\
such that the stabilizer of every line is trivial). A stronger
theorem along the same lines was proved in
\cite[Corollary~52]{LML08}. Furthermore, in~\cite{IK18} the geometry
of van Lint and Schrijver was listed as a counterexample to the
Manickam--Mikl\'{o}s--Singhi conjecture for partial geometries.

In this note we construct another partial geometry with para\-meters
$pg(5,5,2)$, not isomorphic to the partial geometry of van Lint and
Schrij\-ver. The new partial geometry was discovered by prescribing
automorphism groups and performing computer searches, with
techniques similar to the ones used in~\cite{KVK20} for
quasi-symmetric designs. In Section~\ref{sec2} we describe a computer-free
construction of the new $pg(5,5,2)$ by changing some lines of the
geometry of van Lint and Schrij\-ver. Properties of the new partial
geometry are described in Section~\ref{sec3}.

\section{Construction of the new partial geometry}\label{sec2}

A partial geometry $pg(5,5,2)$ has $v=81$ points and as many lines.
We shall denote the geometry of van Lint and Schrijver by
$\G=(\P,\L)$, where the set of lines~$\L$ consists of $6$-element
subsets of the set of points~$\P$. Two constructions of $\G$ are
given in~\cite{vLS81}, the first using cyclotomy in the finite field
$\F_{81}$, and the second using the dual code of the repetition code
in $\F_3^5$. The first construction does not essentially use
multiplication in $\F_{81}$ and we describe it here in purely linear
algebraic terms.

Let $V$ be a four-dimensional vector space over $\F_3$ and
$\{e_1,e_2,e_3,e_4\}$ a basis. Define $e_5=-\sum_{i=1}^4 e_i$ and
$S=\{0,e_1,e_2,e_3,e_4,e_5\}$. Then $\L=\{x+S \mid x\in V\}$ is the
set of lines of $\G=(V,\L)$. In~\cite{vLS81}, the fifth roots of
unity were used as the vectors $e_1,\ldots,e_5$, but the proof that
this is a $pg(5,5,2)$ only relies on the following observation. If a
linear combination of the vectors in $S$ is zero, then it must be of
the form $\alpha\cdot 0 + \beta\cdot(e_1+\ldots+e_5)$. Clearly this
holds starting from any basis.

There are $\left[{4\atop 3}\right]_3=40$ three-dimensional subspaces
of $V$ and they intersect $S$ in $1$, $2$, $3$ or $4$ vectors. Take
a subspace $N_0\le V$ such that $|N_0\cap S|=3$, e.g.\ $N_0=\langle
e_1, e_2, e_3-e_4\rangle$, and denote its cosets by $N_1=e_3+N_0$
and $N_2=e_3+e_4+N_0$. Note that $|N_1\cap S|=3$ and $|N_2\cap
S|=0$, from which we get intersection sizes of $N_0$ with the lines
of~$\G$:
\[
|N_0 \cap (x+S)|=\left\{\begin{array}{ll} 3, & \mbox{for } x\in
N_0\cup N_2,\\
0, & \mbox{for } x\in N_1.\end{array}\right.
\]
(A three-dimensional subspace $N\le V$ such that $|N\cap S|=2$ is a
$2$-ovoid, i.e.\ intersects every line of $\G$ in exactly $2$
points; see~\cite{JT89}.)

We now define a new incidence structure $\G'=(V,\L')$ by changing
the lines of $\G$ disjoint from $N_0$, and retaining the $3$-secants
of $N_0$. Take $S'=\{0,e_1',e_2',e_3',e_4', e_5\}$ for $e_1'=
-e_1+e_3$, $e_2'=-e_1+e_3-e_4$, $e_3'= -e_2+e_4$,
$e_4'=-e_2-e_3+e_4$ and define
\[
\L'=\{x+S'\mid x\in N_1\} \cup \{x+S\mid x\in N_0\cup N_2\}.
\]
The set $S'$ is of the same form as $S$, i.e.\
$\{e_1',e_2',e_3',e_4'\}$ is a basis of~$V$ and $e_5=-\sum_{i=1}^4
e_i'$. It has the same intersection pattern with the cosets
of~$N_0$, namely $|N_0\cap S'|=|N_1\cap S'|=3$, $|N_2\cap S'|=0$.
From this and the observation $\Delta(S) \cap \Delta(S')\cap
N_0=\emptyset$ it follows that $\G'$ is a partial linear space with
points and lines of degree $6$. Here $\Delta(S)=\{x-y\mid x,y\in S,
x\neq y\}$ is the set of differences of~$S$. To prove that $\G'$ is
a partial geometry with $\alpha=2$, it suffices to check that the
point graph $\Gamma_1(\G')$ is strongly regular with parameters
$srg(81,30,9,12)$ by the Lemma from~\cite{AMC81}. Because
$\Delta(S) \cap (N_1 \cup N_2) = \Delta(S') \cap (N_1 \cup N_2)$ holds,
collinearity in $\G$ and $\G'$ differs only for points $x,y\in N_1$
and $x,y\in N_2$. In all other cases the common neighbours
of $x$ and $y$ are the same as in $\Gamma_1(\G)$, which we know is a
$srg(81,30,9,12)$. If $x,y\in N_1$ are collinear, one
can check that they have three common neighbours in $N_1$.
Outside of $N_1$ the common neighbours are the same as in
$\Gamma_1(\G)$, hence there are $6$ of them. It follows that
two adjacent vertices of $\Gamma_1(\G')$ always have $\lambda=9$
common neighbours. With a little more effort one can check that
two non-adjacent vertices $x,y\in N_1$ always have $\mu=12$ common
neighbours, and the arguments are the same for $x,y\in N_2$.

Finally, to prove that the partial geometries $\G$ and $\G'$ are not
isomorphic, take two collinear points $x=0$ and $y=e_1$. The set of
points collinear with both $x$ and $y$ is $A\cup B\cup \{z\}$ for
$A=\{e_2,e_3,e_4,e_5\}$,
$B=\{e_1-e_2,e_1-e_3,e_1-e_4,-e_1+e_2+e_3+e_4\}$, and $z=-e_1$ (in
both geometries). In~\cite{CvL82} it was noted that the collinearity
graph of this subconfiguration of $\G$ always consists of two copies
of the complete graph $K_4$ and an isolated vertex. Indeed, pairs of
points in $A$ and pairs of points in $B$ are collinear and there are
no other collinearities between points of $A\cup B\cup \{z\}$
in~$\G$. In $\G'$ we lose mutual collinearity of $e_1-e_3$,
$e_1-e_4$, and $-e_1+e_2+e_3+e_4$ and gain no further
collinearities, so $B$ is a star graph.

\section{Properties of the new partial geometry}\label{sec3}

In~\cite{CvL82}, the full automorphism group of $\G$ was determined
as $\Aut(\G)=\F_3^4\rtimes S_6$ (semidirect product with normal
subgroup $\F_3^4$) of order $58\,320$. The new geometry $\G'$
clearly has $N_0=\F_3^3$ as automorphism group. Using
nauty~\cite{MP14} and GAP~\cite{GAP4}, we found that the full
automorphism group is $\Aut(\G')=\F_3^3\rtimes G$ of order $972$,
where $G\le S_6$ is a subgroup of order $36$ isomorphic to
$\F_3^2\rtimes \Z_4$. The group $\Aut(\G')$ is not transitive and
does not have Singer subgroups. The point orbits of $\Aut(\G')$ are
$N_0$ and $N_1\cup N_2$, and the line orbits are $\{x+S'\mid x\in
N_1\}$ and $\{x+S\mid x\in N_0\cup N_2\}$.

The point graphs $\Gamma_1(\G)$ and $\Gamma_1(\G')$ are strongly
regular with parameters $srg(81,30,9,12)$ and are not isomorphic.
The full automorphism group $\Aut \Gamma_1(\G)$ is twice as large as
$\Aut(\G)$~\cite{CvL82}, while the full automorphism groups $\Aut
\Gamma_1(\G')$ and $\Aut(\G')$ coincide. Using Cliquer~\cite{NO03},
we found that the graphs $\Gamma_1(\G)$ and $\Gamma_1(\G')$ are
faithfully geometric, meaning that they only support their
respective partial geometries up to isomorphism. The graph
$\Gamma_1(\G)$ has $162$ cliques of size $6$ and the graph
$\Gamma_1(\G')$ has $108$ such cliques.

The two geometries are self-dual, hence their line graphs
$\Gamma_2(\G)$ and $\Gamma_2(\G')$ are isomorphic to the point
graphs. Of the $162$ cliques of size $6$ in $\Gamma_2(\G)$, $81$
correspond to stars (sets of $6$ lines through a single point). The
remaining $81$ cliques correspond to negative lines, i.e.\ sets of
the form $-(x+S)$, $x\in V$. For every negative line there are $6$
lines of $\G$ intersecting it in precisely one point and they form a
clique in $\Gamma_2(\G)$. Negative lines were used in~\cite{IK18} to
show that $\G$ does not have the strict MMS star property. There are
$27$ cliques of size $6$ in the line graph $\Gamma_2(\G')$ that are
not stars. They can be used in the same way to show that the new
partial geometry $\G'$ is also a counterexample to the MMS
conjecture.

\section*{Note added in proof}

After submitting the paper, the author learned that the new partial
geometry was also independently discovered in~\cite{CST21} by
a different method.

\end{document}